\documentstyle[a4]{article}%
\input amssym.def
\input amssym
\font\teneusm=eusm10
\font\seveneusm=eusm7
\font\fiveeusm=eusm5
\newfam\eusmfam
\textfont\eusmfam=\teneusm
\scriptfont\eusmfam=\seveneusm
\scriptscriptfont\eusmfam=\fiveeusm
\def\script#1{{\fam\eusmfam\relax#1}}
\def\Z{\Bbb Z}
\def\T{\Bbb T}
\def\R{\Bbb R}

\def\C{\Bbb C}

\def\e{\mathop{\fam0 e}\nolimits}
\def\ii{{\fam0 i}} 
\def\se{\subseteq}
\def\dem{{\it Proof.\/ }} 
\def\bull{\leavevmode\kern .1ex
\vrule height 1.2ex width 1.1ex depth-.1ex \kern.1ex} 
\def\Ref#1{\hbox{$(\ref{#1})$}}
\def\({{\fam0\rm (}}
\def\){\/{\fam0\rm )}} 
\def\ie{{\it i.\,e.}\ } 
\def\vs{{\it vs.}\ } 
\def\th{\vartheta}
\def\rh{\varrho}

\newtheorem{thm}{Theorem}[section]
\newtheorem{lem}[thm]{Lemma}

\newtheorem{prp}[thm]{Proposition}
\newtheorem{cor}[thm]{Corollary}
\newtheorem{dfn}[thm]{Definition}

\def\exa{
\addtocounter{thm}{1}{\bf Example \thethm}\hskip1ex}

\title{The Sidon constant of sets with three elements}
\author{Stefan Neuwirth\footnotemark}\date{}

\begin{document}
\maketitle

\footnotetext[1]{
\noindent Laboratoire de Math\'ematiques, Universit\'e de
Franche-Comt\'e, 25030 Besan\c con cedex
}
\begin{abstract}
\noindent 
We solve an elementary minimax problem and obtain the exact value of
the Sidon constant for sets with three elements $\{n_0,n_1,n_2\}$: it is 
$\sec(\pi/2n)$ for 
$\displaystyle n={\max|n_i-n_j|/\gcd(n_1-n_0,n_2-n_0)}$.
\end{abstract}
\section{Introduction}
Let $\Lambda=\{\lambda_0,\lambda_1,\lambda_2\}$ be a set of three
frequencies and $\varrho_0,\varrho_1,\varrho_2$ three
positive intensities. We solve the following extremal problem: 

\vskip\abovedisplayskip
$(\dagger)\quad\vcenter{
\halign{#\hfil\cr
To find $\vartheta_0,\vartheta_1,\vartheta_2$ three phases such that,
putting $c_j=\varrho_j\e^{\ii\vartheta_j}$,\cr 
the maximum $\max_t|c_0\e^{\ii \lambda_0t} + c_1\e^{\ii \lambda_1t} +
c_2\e^{\ii \lambda_2t}|$ is minimal.\cr}
}$
\vskip\belowdisplayskip

This enables us to generalize a result of D.\ J.\ Newman. He solved the
following extremal problem for $\Lambda=\{0,1,2\}$:

\vskip\belowdisplayskip$(\ddagger)\quad\vcenter{\halign{#\hfil\cr
To find $f(t)=c_0\e^{\ii \lambda_0t}
+c_1\e^{\ii \lambda_1t} +c_2\e^{\ii \lambda_2t}$ with
$\|f\|_\infty=\max_t|f(t)|\le1$\cr 
such that $\|\widehat{f}\|_1=|c_0|+|c_1|+|c_2|$ is maximal.\cr}}$
\vskip\belowdisplayskip 

Note that for such an
$f$, $\|\widehat{f}\|_1$ is the Sidon constant of $\Lambda$. Newman's
argument is the following (see \cite[Chapter 3]{sh51}): by the
parallelogram law, 
\begin{eqnarray*}
\max_t|f(t)|^2&=&\max_t|f(t)|^2\vee|f(t+\pi)|^2\\
&\ge&\max_t\bigl(|f(t)|^2+|f(t+\pi)|^2\bigr)/2\\
&=&\max_t\bigl(|c_0+c_1\e^{\ii t}+c_2\e^{\ii 2t}|^2+
               |c_0-c_1\e^{\ii t}+c_2\e^{\ii 2t}|^2\bigr)/2\\
&=&\max_t|c_0+c_2\e^{\ii2t}|^2+|c_1|^2=
\bigl(|c_0|+|c_2|\bigr)^2+|c_1|^2\\
&\ge&\bigl(|c_0|+|c_1|+|c_2|\bigr)^2/2
\end{eqnarray*}
and equality holds exactly for multiples and translates of $f(t)=1+2\ii\e^{\ii
t}+\e^{\ii2t}$.

Let us describe this paper briefly. We use a real-variable approach:
Problem $(\dagger)$ reduces to studying a function of form 
$$
\Phi(t,\vartheta)=|1+r\e^{\ii\vartheta}\e^{\ii kt}+s\e^{\ii
lt}|^2\hbox{ for }r,s>0,\ k\ne l\in\Z^*
$$
and more precisely $\Phi^*(\vartheta)=\max_t\Phi(t,\vartheta)$. We
obtain the variations of $\Phi^*$: the point is that we find ``by
hand'' a local minimum of $\Phi^*$ and that any two minima of $\Phi^*$
are separated by a maximum of $\Phi^*$, which corresponds to an
extremal point of $\Phi$ and therefore has a handy description. The
solution to Problem $(\ddagger)$ then turns out to derive easily from
this. 

The initial motivation was twofold. In the first place, we wanted to
decide whether sets 
$\Lambda=\{\lambda_n\}$ such that $\lambda_{n+1}/\lambda_n$ is bounded
by some $q$ may have a Sidon constant arbitrarily close to $1$ and to
find evidence among sets with three elements. That there are such
sets, finite but arbitrarily large, may in fact 
be proven by the method of Riesz products in \cite[Appendix V,
\S1.II]{ks63}. In the second place, we wished to show that 
the real and complex unconditionality constants are distinct for basic
sequences of characters $\e^{\ii nt}$; we prove however that they
coincide in the space $\script{C}(\T)$ for sequences with three terms. 

\vskip\belowdisplayskip
{\bf Notation } $\T=\{z\in\C:|z|=1\}$ and $\e_\lambda(z)=z^\lambda$
for $z\in\T$ and $\lambda\in\Z$. 
\section{Definitions}
\begin{dfn}
$(1)$ Let $\Lambda\se\Z$. $\Lambda$ is a Sidon set if there is a
constant $C$ such that for all trigonometric polynomials
$f(t)=\sum_{\lambda\in\Lambda}c_\lambda\e^{\ii\lambda t}$ with
spectrum in $\Lambda$ we have 
$$
\|\widehat{f}\|_1=\sum_{\lambda\in\Lambda}|c_\lambda|\le
C\max_t|f(t)|=\|f\|_\infty.
$$
The optimal $C$ is called the Sidon constant of $\Lambda$.

$(2)$ Let $X$ be a Banach space. The sequence $(x_n)\se X$ is a real
\(\vs complex\) unconditional basic sequence in $X$ if there is a
constant $C$ 
such that 
$$\Bigl\|\sum\vartheta_nc_nx_n\Bigr\|_X\le C\Bigl\|\sum c_nx_n\Bigr\|_X$$
for every real \(\vs complex\) choice of signs $\vartheta_n\in\{-1,1\}$
\(\vs $\vartheta_n\in\T$\)  
and every finitely supported family of coefficients $(c_n)$. The
optimal $C$ is the real \(\vs complex\) unconditionality constant of
$(x_n)$ in $X$.
\end{dfn}
Let us state the two following well known facts.
\begin{prp}
$(1)$ The Sidon constant of $\Lambda$ is the complex unconditionality
constant of the sequence of functions 
$(\e_\lambda)_{\lambda\in\Lambda}$ in the space $\script{C}(\T)$.

$(2)$ The complex unconditionality constant is at most $\pi/2$ times
the real unconditionality constant.
\end{prp}
\dem
$(1)$ holds because
$\bigl\|\sum\vartheta_\lambda c_\lambda\e_\lambda\bigr\|_\infty=
\sum|c_\lambda|$ for
$\vartheta_\lambda=\overline{c_\lambda}/|c_\lambda|$. 

$(2)$ Because the complex unconditionality constant of the sequence
$(\epsilon_n)$ of
Rade\-macher functions in $\script{C}(\{-1,1\}^\infty)$ is $\pi/2$
(see \cite{se97}), 
\begin{eqnarray*}
\sup_{\vartheta_n\in\T}\Bigl\|\sum\vartheta_nc_nx_n\Bigr\|_X
&=&\sup_{x^*\in B_{X^*}}\sup_{\vartheta_n\in\T}\sup_{\epsilon_n=\pm1}
\Bigl|\sum\vartheta_n c_n\langle x^*,x_n\rangle\epsilon_n\Bigr|\\&\le&
\pi/2\sup_{x^*\in  B_{X^*}}\sup_{\epsilon_n=\pm1}
\Bigl|\sum c_n\langle x^*,x_n\rangle\epsilon_n\Bigr|\\
&=&\pi/2\sup_{\epsilon_n=\pm1}\Bigl\|\sum \epsilon_nc_nx_n\Bigr\|_X.
\end{eqnarray*}
Furthermore the real unconditionality constant of $(\epsilon_n)$ in
$\script{C}(\{-1,1\}^\infty)$ is $1$: therefore the factor $\pi/2$ is
optimal.\nolinebreak\hfill\bull\vskip\belowdisplayskip 

Let us straighten out the expression of the Sidon constant. For 
$$f(t)=c_0\e^{\ii \lambda_0t}+c_1\e^{\ii \lambda_1t}+
c_2\e^{\ii\lambda_2t}\ ,\ c_j=\varrho_j\e^{\ii\vartheta_j}\ ,$$
the supremum norm $\|f\|_\infty$ of $f$ is equal to 
\begin{equation}
\label{ubc:angle}
\|\varrho_0+\varrho_1\e^{\ii\vartheta}\e_{\lambda_1-\lambda_0}+
\varrho_2\e_{\lambda_2-\lambda_0}\|_\infty\ ,\ 
\vartheta={\lambda_1-\lambda_2\over\lambda_2-\lambda_0}\vartheta_0
+\vartheta_1+{\lambda_0-\lambda_1\over\lambda_2-\lambda_0}\vartheta_2
\end{equation}
and therefore the Sidon constant $C$ of
$\Lambda=\{\lambda_0,\lambda_1,\lambda_2\}$ may be written
\begin{equation}
  \label{ubc:ubc}
C=\max_{r,s>0,\vartheta\in\T}
(1+r+s)/\|1+r\e^{\ii\vartheta}\e_k+s\e_l\|_\infty\quad\hbox{with }
\left\{
\begin{array}{l}
k=\lambda_1-\lambda_0\\
l=\lambda_2-\lambda_0.
\end{array}\right.  
\end{equation}
By change of variables, we may suppose w.l.o.g.\ that $k$ and $l$ are
coprime. 
\section{A solution to Extremal problem $(\dagger)$}
Let us first establish
\begin{lem}\label{ubc:lem}
Let $\lambda_1,\dots,\lambda_k\in\Z^*$ and $\rh_1,\dots,\rh_k>0$.
Let 
$$
f(t,\th)=1+\rh_1\e^{\ii(\lambda_1t+\th_1)} + \dots +
\rh_{k-1}\e^{\ii(\lambda_{k-1}t+\th_{k-1})} +
\rh_k\e^{\ii\lambda_kt}
$$
and $\Phi(t,\th)=|f(t,\th)|^2$. The extremal points $(t,\th)$
such that $\nabla\Phi(t,\vartheta)=0$ satisfy either
$f(t,\th)=0$ or
$\lambda_1t+\th_1 \equiv \dots \equiv \lambda_{k-1}t+\th_{k-1} \equiv
\lambda_kt \equiv0\ \hbox{mod.}\ \pi$.
\end{lem}
\dem
As $\Phi=(\Re f)^2+(\Im f)^2$, the extremal points
$(t,\th)$ satisfy 
$$
\left\{\vcenter{\halign{\hfil$#$&\,\hfil$#$\cr
\Re{\partial f\over\partial t}(t,\th)\,\Re f(t,\th)\,+&
\Im{\partial f\over\partial t}(t,\th)\,\Im f(t,\th)=0\cr
-\sin(\lambda_i t+\th_i)\,\Re f(t,\th)\,+&
\cos(\lambda_i t+\th_i)\,\Im f(t,\th)=0\cr}}\right.
$$
for $1\le i\le k-1$. Suppose that $f(t,\th)\ne0$: then the
system above implies that $\lambda_1t+\th_1 \equiv \dots \equiv
\lambda_{k-1}t+\th_{k-1}\equiv \tau\ \hbox{mod.}\ \pi$ for a certain
$\tau$, so that this system becomes

$$\left\{\vcenter{\halign{\hfil$#$&\,\hfil$#$\cr
(-\sigma\sin\tau-\lambda_k\rh_k\sin\lambda_kt)\,\Re f(t,\th)\,+& 
( \sigma\cos\tau+\lambda_k\rh_k\cos\lambda_kt)\,\Im f(t,\th)=0\cr
-\sin\tau\,\Re f(t,\th)\,+&\cos\tau\,\Im f(t,\th)=0\cr}}\right.$$

with
$\sigma=\lambda_1\rh_1+\dots+\lambda_{k-1}\rh_{k-1}$. Therefore
also $\lambda_kt\equiv0\ \hbox{mod.}\ \pi$, so that 
$$
-\sin\tau\cdot(1+(\rh_1+\dots+\rh_k)\cos\tau)+\cos\tau\cdot(\rh_1+\dots+\rh_k)\sin\tau=0
$$
and $\tau\equiv 0\ \hbox{mod.}\ \pi$.\hfill\bull\vskip\belowdisplayskip

The following result is the core of the paper.
\begin{lem}\label{ubc:prp}
Let $r,s>0$, $k,l\in\Z^*$ distinct and coprime. Let
\begin{eqnarray*}
\Phi(t,\vartheta)&=&|1+r\e^{\ii\vartheta}\e^{\ii kt}+s\e^{\ii lt}|^2\\
&=&1+r^2+s^2+2r\cos(kt+\vartheta)+2s\cos{lt}+2rs\cos((l-k)t-\vartheta).
\end{eqnarray*}
Let $\Phi^*(\vartheta)=\max_t\Phi(t,\vartheta)$. Then
$\Phi^*$ is even, has period $2\pi/|l|$ and decreases on
$[0,\pi/|l|]$. Therefore $\min_\vartheta\Phi^*(\vartheta)=\Phi^*(\pi/l)$. 
\end{lem}
\dem
%Let us first locate 
By Lemma \ref{ubc:lem}, the extremal points of
$\Phi$ satisfy either $\Phi(t,\vartheta)=0$ or 
%. $\nabla\Phi(t,\vartheta)=0$ is Equation 
%\Ref{ubs:S} with 
%$$(a_{ij})=\pmatrix{
%rk&sl&rs(l-k)\cr
%r&0&-rs\cr}\quad\hbox{and}\quad\cases{\alpha=kt+\vartheta&\cr\beta=lt.\cr}
%$$ 
%In the notation of Lemma \ref{ubc:lem}, $d_1=-rsl\ne0$, $d_2=-r^2sl\ne0$,
%$d_3=-rs^2l\ne0$. The solution \Ref{ubc:sol}, if it exists, yields
%$\Phi(t,\vartheta)=0$ and corresponds to the absolute minimum of
%$\Phi$. Every other extremal point of $\Phi$ satisfies 
$$
  \label{ubc:lem:cons}
kt+\vartheta\equiv lt\equiv0\hbox{ mod.\ }\pi. 
$$

$\Phi^*$ is continuous (see \cite[Chapter 5.4]{ps72}) and
$(2\pi/|l|)$-periodical: choose $j\in\Z$ such that $jk\equiv1$ mod.\
$l$. Then
$$\Phi(t+2j\pi/l,\vartheta) =
|1+r\e^{\ii(\vartheta+2\pi jk/l)}\e^{\ii kt}+s\e^{\ii lt}|^2 = 
\Phi(t,\vartheta+2\pi/l).$$ 
Furthermore 
$\Phi(t,-\vartheta)=\Phi(-t,\vartheta)$ and $\Phi^*$ is even. Thus
$\Phi^*$ attains its minimum on $[0,\pi/|l|]$.

Let us show that $\Phi^*$ has a local minimum at $\pi/l$ for all
values of $r,s>0$ except eventually one value of $s$ for a given
$r$. Let $t^*$ be such that $\Phi^*(\pi/l)=\Phi(t^*,\pi/l)$. As
before, for $j$ such that $jk\equiv-1$ mod.\ $l$,
$$\Phi(t^*,\pi/l+\vartheta) = \Phi(-t^*,-\pi/l-\vartheta) =
\Phi(-t^*+2j\pi/l,\pi/l-\vartheta).$$ 
If $\partial\Phi/\partial\vartheta(t^*,\pi/l)\ne0$, this shows that
$\Phi^*$ has a local minimum and a cusp at $\pi/l$. Let us suppose that
$\partial\Phi/\partial\vartheta(t^*,\pi/l)=0$. Then 
$\nabla\Phi(t^*,\pi/l)=0$ and 
therefore $kt^*+\pi/l=j\pi$ and $lt^*=j'\pi$ for some
$j,j'\in\Z$. Then $j$ or $j'$ must be odd. 
We have 
\begin{equation}\label{ubc:d2phi}
{\partial^2\Phi\over\partial t^2}(t^*,\pi/l) =
-2\bigl(rk^2(-1)^j+sl^2(-1)^{j'}+rs(l-k)^2(-1)^{j+j'}\bigr).
\end{equation}
We shall suppose that $\partial^2\Phi/\partial t^2(t^*,\pi/l)\ne0$,
which removes at most one value of $s$ for a given $r$. Then, by
the Theorem of Implicit Functions, there is a unique differentiable
function $t_*$ defined  in a neighbourhood of $\pi/l$ such
that 
$$t_*(\pi/l)=t^*\ ,\ \partial\Phi/\partial
t(t_*(\vartheta),\vartheta)=0\ ,\ \partial^2\Phi/\partial
t^2(t_*(\vartheta),\vartheta)<0.$$  
Let
$\Phi_*(\vartheta)=\Phi(t_*(\vartheta),\vartheta)$. Then we have
$\Phi_*'(\pi/l)=\partial\Phi/\partial\vartheta(t^*,\pi/l)=0$ and a
computation yields 
$$\Phi_*''(\pi/l) = { 
{\partial^2\Phi\over\partial t^2} {\partial^2\Phi\over\partial\vartheta^2} -
\bigl({\partial^2\Phi\over\partial t\partial\vartheta}\bigr)^2 \over 
{\partial^2\Phi\over\partial t^2}} (t^*,\pi/l) = 4rsl^2
\biggl({\partial^2\Phi\over\partial t^2}(t^*,\pi/l)\biggr)^{-1}\Delta,$$
where $\Delta=(-1)^{j+j'}+r(-1)^{j'}+s(-1)^j$.
Let us prove that $\Delta<0$ and thus
$\Phi_*''(\pi/l)>0$ and that therefore $\Phi_*$ and consequently
$\Phi^*$ have a local minimum at 
$\pi/l$. If we had $\Delta\ge0$ and

\bull\  $j$ even, $j'$ odd: then $-1-r+s\ge0$ and by \Ref{ubc:d2phi}
$${\partial^2\Phi\over\partial t^2}(t^*,\pi/l) \ge
2\bigl(-rk^2+(1+r)l^2+r(1+r)(l-k)^2\bigr)=2(r(l-k)+l)^2\ge0;$$ 

\bull\  $j$ odd, $j'$ even: then $-1+r-s\ge0$ and by \Ref{ubc:d2phi}
$${\partial^2\Phi\over\partial t^2}(t^*,\pi/l) \ge
2\bigl((1+s)k^2-sl^2+s(1+s)(l-k)^2\bigr)=2(s(l-k)-k)^2\ge0;$$  

\bull\  $j$ odd, $j'$ odd: then $1-r-s\ge0$. Considering
\Ref{ubc:d2phi}, we have the following alternative. 
If $l^2\ge r(l-k)^2$, then 
$\partial^2\Phi/\partial t^2(t^*,\pi/l) \ge0$; otherwise
$${\partial^2\Phi\over\partial t^2}(t^*,\pi/l) \ge 
2\bigl(rk^2+(1-r)(l^2-r(l-k)^2)\bigr)=2(r(l-k)-l)^2\ge0.$$
Let us show that then $\Phi^*$ must decrease on $[0,\pi/|l|]$. Otherwise
there are $0\le\vartheta_0<\vartheta_1\le\pi/|l|$ such that
$\Phi^*(\vartheta_1)>\Phi^*(\vartheta_0)$. As $\pi/|l|$ is a local
minimum, there is a $\vartheta_0<\vartheta^*<\pi/|l|$ such that 
$$\Phi^*(\vartheta^*) =
\max_{\vartheta_0\le\vartheta\le\pi/|l|}\Phi^*(\vartheta) =
\max_{\scriptstyle          0\le t        < 2\pi\atop
      \scriptstyle\vartheta_0\le\vartheta\le \pi/|l|}
\Phi(t,\vartheta),$$
\ie there further is some $t^*$ such that $\Phi$ has a local maximum
at $(t^*,\vartheta^*)$. But then 
$kt^*+\vartheta^*\equiv lt^*\equiv0$ mod.\ $\pi$
and $\vartheta^*\equiv0$ mod.\ $\pi/l$ and this is
false. That shows the proposition, except for one
value of $s$ at most for a given $r$. But $\Phi^*$ is a continuous
function of $s$ and the proposition is true by a
perturbation.\nolinebreak\hfill\bull\vskip\belowdisplayskip
By Computation \Ref{ubc:angle} and Lemma \ref{ubc:prp}, we obtain 
\begin{thm}
Let $\lambda_0,\lambda_1,\lambda_2\in\R$ and
$\varrho_0,\varrho_1,\varrho_2>0$. The solution to Extremal
problem $(\dagger)$ is the following.

\bull\  If the smallest additive group containing
$\lambda_1-\lambda_0$ and $\lambda_2-\lambda_0$ is dense in $\R$,
then the maximum is independent of the phases
$\vartheta_0,\vartheta_1,\vartheta_2$ and makes
$\varrho_0+\varrho_1+\varrho_2$.  

\bull\  Otherwise let $d=\gcd(\lambda_1-\lambda_0,\lambda_2-\lambda_0)$
be a generator of this group. Then the sought phases
$\vartheta_0,\vartheta_1,\vartheta_2$ are given by
$$
\vartheta_0(\lambda_2-\lambda_1)+
\vartheta_1(\lambda_0-\lambda_2)+
\vartheta_2(\lambda_1-\lambda_0)
\equiv d\pi\hbox{ mod.\ }2d\pi.
$$
In particular, these phases may be chosen among $0$ and $\pi$.
\end{thm}
\section{A solution to Extremal problem $(\ddagger)$}
There are two cases where one can make explicit computations by Lemma
\ref{ubc:prp}.  

\vskip\abovedisplayskip\exa
The real and complex unconditionality constant of
$\{0,1,2\}$ in $\script{C}(\T)$  is $\sqrt{2}$. Indeed, a case study shows that
$$\|1+\ii r\e_1+s\e_2\|_\infty = 
\cases{r+|s-1|&if $r|s-1|\ge4s$\cr
(1+s)(1+r^2/4s)^{1/2}&if $r|s-1|\le4s$\cr}$$
and this permits to compute the maximum \Ref{ubc:ubc}, which is obtained
for $r=2$, $s=1$. This yields another proof to Newman's result
presented in the Introduction.\vskip\belowdisplayskip

\exa
The real and complex unconditionality constant of
$\{0,1,3\}$ in $\script{C}(\T)$ is $2/\sqrt{3}$. Indeed, a case study shows that
$\|1+r\e^{\ii \pi/3}\e_1+s\e_3\|_\infty$ makes 
$$\cases{
1+r-s&if $s\le r/(4r+9)$\cr
\bigl({2\over27}s(r^2+9+3r/s)^{3/2}
-{2\over27}r^3s+{2\over3}r^2+rs+s^2+1\bigr)^{1/2}
&if $s\ge r/(4r+9)$\cr}$$
and this permits to compute the maximum \Ref{ubc:ubc}, which is obtained
exactly at  $r=3/2$, $s=1/2$.\vskip\belowdisplayskip

These examples are particular cases of the following theorem.

\begin{thm}\label{ubc:THM}
Let $\lambda_0,\lambda_1,\lambda_2\in\Z$ be distinct. Then the 
Sidon\index{Sidon set!constant} constant of
$\Lambda=\{\lambda_0,\lambda_1,\lambda_2\}$ is 
$\sec(\pi/2n)$, where
$n=\max|\lambda_i-\lambda_j|/\gcd(\lambda_1-\lambda_0,\lambda_2-\lambda_0)$.

\end{thm}
\dem
We may suppose $\lambda_0<\lambda_1<\lambda_2$. Let $k=(\lambda_1-\lambda_0)/\gcd(\lambda_1-\lambda_0,\lambda_2-\lambda_0)$ and $l=(\lambda_2-\lambda_0)/\gcd(\lambda_1-\lambda_0,\lambda_2-\lambda_0)$. 
By Lemma \ref{ubc:prp},
the Arithmetic-Geometric Mean Inequality bounds the Sidon constant $C$ of
$\{0,k,l\}$ in the following way:
\begin{eqnarray*}
C=\max_{r,s>0}{1+r+s\over\|1+r\e^{\ii\pi/l}\e_k+s\e_l\|_\infty}&\le&
\max_{r,s>0}{1+r+s\over|1+r\e^{\ii\pi/l}+s|}\\&=&
\max_{r,s>0}\biggl(1-\sin^2{\pi\over2l}{4r(1+s)\over(1+r+s)^2}\biggr)^{-1/2}\\
&\le&
\bigl(1-\sin^2(\pi/2l)\bigr)^{-1/2}=\sec(\pi/2l).
\end{eqnarray*}
This inequality is sharp: we have equality for $s=k/(l-k)$ and
$r=1+s$. In fact the derivative of 
$|1+r\e^{\ii\pi/l}\e^{\ii kt}+s\e^{\ii lt}|^2$ is
then
$$
{8kl\over k-l}\cos{kt+\pi/l\over2}\sin{lt\over2}
\cos{(l-k)t-\pi/l\over2}
$$
and its extremal points are 
$$
{2j+1\over   k}\pi-{\pi\over kl}\ ,\ 
{2j  \over l  }\pi\ ,\ 
{2j+1\over l-k}\pi+{\pi\over l(l-k)}\ :\ j\in\Z$$
so that its extremal values are
$$4s^2\sin^2{2j+1+l\over2k    }\pi\ ,\ 
  4r^2\cos^2{2j+1  \over2l    }\pi\ ,\ 
  4   \cos^2{2j+1+k\over2(l-k)}\pi\ :\ j\in\Z.$$
Therefore the maximum of $|1+r\e^{\ii\pi/l}\e^{\ii kt}+s\e^{\ii lt}|$ is
$2r\cos(\pi/2l)$. \nolinebreak\hfill\bull\vskip\belowdisplayskip
 
This proof and \Ref{ubc:angle} yield also the
  more precise 
\begin{prp}
Let $\Lambda=\{\lambda_0,\lambda_1,\lambda_2\}\se\Z$. The solution to
Extremal problem
$(\ddagger)$ is a multiple of
%The Sidon constant of $\{n_1,n_1,n_2\}$ is attained for 
$$
f(t)=
\epsilon_0\,|\lambda_1-\lambda_2|\,\e^{\ii \lambda_0t}+
\epsilon_1\,|\lambda_0-\lambda_2|\,\e^{\ii \lambda_1t}+
\epsilon_2\,|\lambda_0-\lambda_1|\,\e^{\ii \lambda_2t}
$$
with $\epsilon_0,\epsilon_1,\epsilon_2\in\{-1,1\}$ real signs such
that

\bull\  
$\epsilon_0\epsilon_1=-1$ if $2^j\mid \lambda_1-\lambda_0$ and
$2^j\nmid \lambda_2-\lambda_0$ for some $j$;

\bull\  
$\epsilon_0\epsilon_2=-1$ if $2^j\nmid \lambda_1-\lambda_0$ and
$2^j\mid \lambda_2-\lambda_0$ for some $j$;

\bull\  
$\epsilon_1\epsilon_2=-1$ otherwise.

The Sidon constant of $\Lambda$ is attained for this $f$. Therefore
the 
complex and real unconditionality constants of
$\{\e_\lambda\}_{\lambda\in\Lambda}$ in $\script{C}(\T)$ 
coincide for sets $\Lambda$  with three elements.
\end{prp}

\section{Some consequences}

Let us underline the following consequences of our computation. 
\begin{cor}
$(1)$ The Sidon constant of sets with three elements is at most
$\sqrt{2}$. 

$(2)$ The Sidon constant of $\{0,n,2n\}$ is $\sqrt{2}$, while the
Sidon constant of $\{0,n+1,2n\}$ is at most
$\sec(\pi/2n)=1+\pi^2/8n^2+o(n^{-2})$ and thus arbitrarily close to $1$.  

$(3)$ The Sidon constant of $\{\lambda_0<\lambda_1<\lambda_2\}$ does
not depend on $\lambda_1$ but on the g.c.d.\ of $\lambda_1-\lambda_0$
and $\lambda_2-\lambda_0$.
\end{cor}

Theorem \ref{ubc:THM} also shows anew that no set of integers with
more than two elements has Sidon constant $1$ (see \cite[p.\ 21]{sh51} or \cite{chm81}). 
Recall now that
$\Lambda=\{\lambda_n\}\se\Z$ is a Hadamard set if there is 
a $q>1$ such that $|\lambda_{n+1}/\lambda_n|\ge q$ for all $n$.
By \cite[Cor.\ 9.4]{ne98}, the Sidon constant of
$\Lambda$ is at most $1+\pi^2/(2q^2-2-\pi^2)$ if 
$q>\sqrt{\pi^2/2+1}\approx2.44$. On the other hand
Theorem \ref{ubc:THM} shows 
\begin{cor}
$(1)$ If there is an integer $q\ge2$ such
that $\Lambda\supseteq\{\lambda,\lambda+\mu,\lambda+q\mu\}$ for some
integers $\lambda$ and $\mu$, then the Sidon constant of $\Lambda$ is
at least 
$$\sec(\pi/2q)>1+\pi^2/(8q^2).$$

$(2)$ In particular, we have the following bounds for the Sidon
constant $C$ of the set $\Lambda=\{q^k\}$, $q\in\Z\setminus\{0,\pm1,\pm2\}$:
$$
1+{\pi^2\over8\max(-q,q+1)^2}< C\le 
1+{\pi^2\over2q^2-2-\pi^2}.$$
\end{cor}
\section{Three questions}
$(1)$ Is there a set $\Lambda$ for which the
real and complex unconditionality constants of
$\{\e_\lambda\}_{\lambda\in\Lambda}$ in $\script{C}(\T)$ differ~? The same
question is open in spaces $L^p(\T)$, $1\le p<\infty$, and even for the
case of three element sets if $p$ is not a small even integer, and
especially for the set $\{0,1,2,3\}$ in any space but
$L^2(\T)$.\vskip\belowdisplayskip  

$(2)$ Let $q>1$. Are there infinite sets $\Lambda=\{\lambda_n\}$ such that
$|\lambda_{n+1}/\lambda_n|\le q$ with Sidon constant
arbitrarily close to $1$~? What about the sequence of integer parts of
the powers of a transcendental number $\sigma>1$~ (see \cite[Cor.\ 2.10,
Prop.\ 3.2]{ne98})~?\vskip\belowdisplayskip

$(3)$ The only set with more than three elements whose Sidon constant
is known is $\{0,1,2,3,4\}$, for which it makes $2$ (see \cite[Chapter
3]{sh51}). Can one compute the Sidon constant of sets with four
elements~? I conjecture that the Sidon constant of $\{0,1,2,3\}$
is $5/3$.  


\begin{thebibliography}{1}

\bibitem{chm81}
D.~I. Cartwright, R.~B. Howlett and J.~R. McMullen, {\em Extreme values for the
  {S}idon constant}, Proc. Amer. Math. Soc. 81 (1981),  531--537.

\bibitem{ks63}
J.-P. Kahane and R.~Salem, {\em Ensembles parfaits et s\'eries
  trigonom\'etriques}, Hermann, 1963.
\newblock Actualit\'es Scientifiques et Industrielles 1301.

\bibitem{ne98}
S.~Neuwirth, {\em Metric unconditionality and {F}ourier analysis}, Studia Math.
  131 (1998),  19--62.

\bibitem{ps72}
G.~P{\'o}lya and G.~Szeg{\H{o}}, {\em Problems and theorems in analysis. {V}ol.
  {I}: {S}eries, integral calculus, theory of functions}, Springer, 1972.
\newblock Grundlehren der mathematischen Wissenschaften 193.

\bibitem{se97}
J.~A. Seigner, {\em Rademacher variables in connection with complex scalars},
  Acta Math. Univ. Comenian. (N.S.) 66 (1997),  329--336.

\bibitem{sh51}
H.~S. Shapiro, {\em Extremal problems for polynomials and power series},
  Master's thesis, Massachusetts Institute of Technology, 1951.

\end{thebibliography}
\end{document}